\documentclass[12pt]{article}

\usepackage[nonatbib, preprint]{neurips_2023} 
\usepackage[sortcites=true]{biblatex}
\usepackage[utf8]{inputenc} 
\usepackage[T1]{fontenc}    
\usepackage{hyperref}       
\usepackage{url}            
\usepackage{booktabs}       
\usepackage{amsfonts}       
\usepackage{nicefrac}       
\usepackage{microtype}      
\usepackage[table]{xcolor}         

\usepackage[utf8]{inputenc} 
\usepackage[T1]{fontenc} 

\usepackage[default]{raleway}  
\usepackage{fontawesome} 

\usepackage{cmbright}
\usepackage[T1]{fontenc}
\usepackage{mathbbol}

\usepackage{amsmath} 
\usepackage{amssymb}
\usepackage{amsthm}

\newcommand{\ie}{\textit{i.e.}\ }
\newcommand{\eg}{\textit{e.g.}\ }
\newcommand{\nb}{\textit{n.b.}\ }
\newcommand{\cf}{\textit{c.f.}\ }

\usepackage{multirow}
\usepackage{colortbl,booktabs} 
\def\ROWCOLOR{black!15!white}
\newcommand{\prox}[1]{\mbox{prox}_{#1}}

\definecolor{BLUE}{rgb}{0.3,0.3,0.9}
\definecolor{RED}{rgb}{0.8,0.05,0.05}
\definecolor{GREEN}{rgb}{0.05,0.5,0.05}
\definecolor{YELLOW}{rgb}{0.9, 0.7, 0.3} 

\usepackage{lipsum}  
\usepackage{fontawesome}
\usepackage{tcolorbox} 

\usepackage{bbm}

\newcommand{\sA}{{\cal A}} 
 
\newcommand{\sC}{{\cal C}}

\newcommand{\sH}{{\cal H}}
\newcommand{\sI}{{\cal I}}

\newcommand{\sM}{{\cal M}}

\newcommand{\sP}{{\cal P}} 
 
\newcommand{\sR}{{\cal R}}

\newcommand{\bbN}{{\mathbb N}} 
 
\newcommand{\bbR}{{\mathbb R}}

\newcommand{\bbP}{{\mathbb P}}

\usepackage{pgfplots}
\usepgfplotslibrary{fillbetween}

\usepackage{amssymb}
\usepackage{multirow}

\newtheorem{remark}{Remark}

\usepackage{tikz}
\usetikzlibrary{calc} 
\usetikzlibrary{arrows.meta}
\usetikzlibrary{fadings}

\usepackage{amsmath}
\DeclareMathOperator*{\argmax}{\mbox{argmax}}
\DeclareMathOperator*{\argmin}{\mbox{argmin}}

\usepackage{fontawesome}
\usepackage{multirow}

\usepackage[Export]{adjustbox}
\usepackage{subfig}
\usetikzlibrary{tikzmark,calc,,arrows,shapes,decorations.pathreplacing}
\tikzset{every picture/.style={remember picture}}
\usetikzlibrary{fit,shapes.misc}
\usetikzlibrary{positioning,backgrounds,spy} 
\usepackage{newfloat}
\usepackage{listings}
\usepackage{wrapfig}
\definecolor{TypalBlueDark}{HTML}{5467b2}
\usepackage{hyperref}
\hypersetup{
    colorlinks=true,       
    linkcolor=TypalBlueDark,          
    citecolor=TypalBlueDark,        
    filecolor=TypalBlueDark,      
    urlcolor=TypalBlueDark           
}

\title{Explainable AI via Learning to Optimize}

\author{%
  Howard Heaton\thanks{The authors contributed equally.} \\
  Typal Academy  
  \And
  Samy Wu Fung\footnotemark[1] \\
  Colorado School of Mines \\
}

\hypersetup{
    pdftitle={Explainable AI via Learning to Optimize},  
    pdfauthor={Howard Heaton},    
    pdfsubject={Explainable AI via Learning to Optimize},
    pdfcreator={Howard Heaton}, 
    pdfproducer={Technology},  
    pdfkeywords={Optimization,
                 Explainable AI,
                 Trustworthy,
                 Certificate,
                 L2O},
}

\begin{document}


\maketitle

\begin{abstract}
Indecipherable black boxes are common in machine learning (ML), but   applications increasingly require explainable artificial intelligence (XAI). The core of XAI is to establish transparent and interpretable data-driven algorithms.
This work provides concrete tools for XAI in situations where prior knowledge must be encoded and untrustworthy inferences flagged. We use the ``learn to optimize'' (L2O) methodology wherein each inference  solves a data-driven optimization problem. Our L2O models are straightforward to implement, directly encode prior knowledge, and yield theoretical guarantees (\eg satisfaction of constraints). We also propose use of interpretable certificates to verify whether model inferences are trustworthy. Numerical examples are provided in the   applications of dictionary-based signal recovery,   CT imaging,  and arbitrage trading of cryptoassets. Code and additional documentation can be found at \href{https://xai-l2o.research.typal.academy}{xai-l2o.research.typal.academy}.
\end{abstract}

\section{Introduction}

\begin{wrapfigure}[14]{r}{0.55\textwidth}
    \centering
    
    \includegraphics[width=0.5\textwidth]{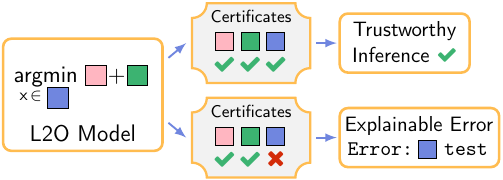}
    \caption{The L2O model is composed of parts (shown as colored blocks) based on prior knowledge or data. L2O inferences solve the optimization problem for given model inputs. Certificates label if each inference is consistent with training. If so, it is trustworthy; otherwise, the faulty model part errs.}
    \label{fig: eye-candy} 
\end{wrapfigure}
A   paradigm shift in machine learning is to construct explainable and transparent models,  often called explainable AI (XAI)  \cite{van2004explainable}. 
This is crucial for sensitive applications like medical imaging and finance (\eg see recent work on the role of explainability  \cite{arrieta2020explainable,adadi2018peeking,dovsilovic2018explainable,samek2019towards}). 
Yet, many commonplace models (\eg fully connected feed forward)   offer limited interpretability. Prior XAI works give explanations via tools like sensitivity analysis \cite{samek2019towards} and layer-wise propagation \cite{bach2015pixel,montavon2019layer}, but these neither quantify trustworthiness nor  necessarily shed light on how to correct ``bad" behaviours. Our work shows how learning to optimize (L2O) can be used to directly embed explainability into models.

The scope of this work is machine learning (ML)   applications where domain experts can create   approximate models by hand. 
In our setting, the inference ${N}_\Theta(d)$ of a model ${N}_\Theta$ with input $d$  solves an optimization problem. That is, we use
\begin{equation}
    {N}_\Theta(d) \triangleq \argmin_{x\in\sC_\Theta(d)} f_\Theta(x; d),
    \label{eq: L2O-model}
\end{equation}
where   $f_\Theta$ is a 
function and $\sC_\Theta(d) \subseteq \bbR^n$ is a constraint set (\eg encoding prior information like   physical quantities), and each  (possibly)   includes dependencies on   weights $\Theta$.
Note the model ${N}_\Theta$ is \textit{implicit} since its output is defined by an optimality condition rather than an explicit computation.
To clarify the scope of the word \textit{explainable} in thi swork, we adopt the following conventions. We say a model is explainable provided a domain expert can identify the core design elements of a model and how they translated to expected inference properties. We say a \textit{particular inference} is explainable provided its properties can be linked to the model's design and intended use. Explainable models and inferences are achieved via L2O with our proposed certificates.

A standard practice in software engineering is to code post-conditions after function calls return. Post-conditions are criteria used to validate what the user expects from the code and ensure   code is not executed under the wrong assumptions. \cite{anaya2018clean} We propose use of these for ML model inferences (see Figures \ref{fig: eye-candy} and \ref{fig: python-sample-code}). 
These conditions enable use of certificates with labels -- pass, warning or fail --  to describe each model inference. We define an inference to be \textit{trustworthy provided it satisfies all provided post-conditions}.

\begin{figure*}[t]
    \centering
    \includegraphics[width=0.95\textwidth]{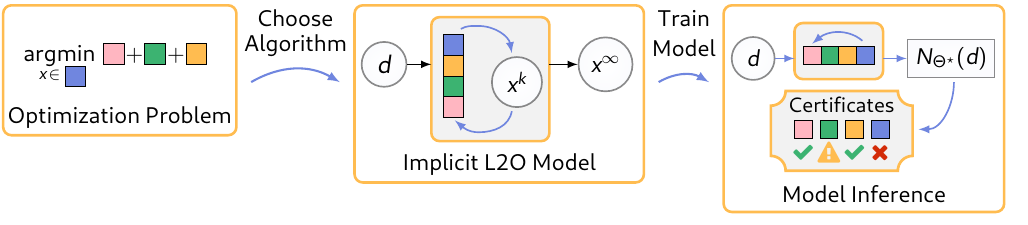}
    \caption{Left shows learning to optimize (L2O) model. Colored blocks denote prior knowledge and data-driven terms. Middle shows an iterative algorithm   formed from the blocks (\eg via proximal/gradient operators) to solve   optimization problem.
    Right   shows a trained model's inference $N_{\Theta^\star}(d)$ and its certificates. Certificates identify if properties of inferences are consistent with training data. 
    Each label is associated with properties of specific blocks (indicated by labels next to blocks in right schematic). Labels take value pass \textcolor{GREEN}{\faCheck}, warning \textcolor{YELLOW}{\faExclamationTriangle}, or fail \textcolor{RED}{\faTimes}, and values  identify if inference features for model parts are trustworthy.}
    \label{fig: certificate-diagram}
\end{figure*} 

\begin{figure*}[t]
  \centering
  \small
     \includegraphics[]{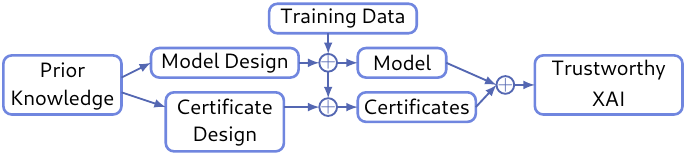}
    \vspace{0.1in}
    \caption{This diagram illustrates   relationships between certificates, models, training data, and prior knowledge. Prior knowledge is embedded directly into model design via the L2O methodology. This also gives rise to quantities to measure for certificate design.  The designed model is   tuned using training data to obtain the ``optimal'' L2O model (shown by arrows touching top middle $+$ sign). The certificates are tuned to match the test samples and/or model inferences on training data (shown by arrows with bottom middle $+$ sign). Together the model and certificates yield inferences with certificates of trustworthiness.}
    \label{fig: diagram-certificate-model-flow}
\end{figure*}

Two ideas, optimization and certificates, form a concrete notion of XAI. 
Prior and data-driven knowledge can be encoded via optimization, and this encoding can be verified via certificates (see Figure \ref{fig: diagram-certificate-model-flow}).   
To illustrate, consider inquiring why a   model generated a ``bad'' inference (\eg an inference disagrees with observed measurements). The first diagnostic step is to check   certificates. If no fails occurred, the model was not designed to handle the instance encountered. In this case, the model in (\ref{eq: L2O-model}) can be redesigned to encode  prior knowledge of the situation. Alternatively, each failed certificate shows a type of error {and}   often corresponds to portions of the model (see Figures \ref{fig: eye-candy} and \ref{fig: certificate-diagram}). The L2O model allows   debugging of algorithmic implementations and assumptions to correct   errors. In a sense, this setup enables one to manually backpropagate errors to fix models (similar to training).

\paragraph{Contributions}
This work brings new explainability and guarantees to deep learning applications   using prior knowledge. We propose novel implicit L2O models with intuitive design, memory efficient training, inferences that satisfy optimality/constraint conditions, and certificates that either indicate trustworthiness or flag inconsistent inference features.  

\begin{table}[t]
    \centering 
    \renewcommand{\arraystretch}{1.25}
    \begin{tabular}{c|c|c|c} 
        L2O & Implicit & Flags & {Obtainable} Model Property \\\hline 
        \rowcolor{\ROWCOLOR}
        \checkmark &  &  & Intuitive Design\\
         & \checkmark &   & Memory Efficient  \\
        \rowcolor{\ROWCOLOR} 
       \checkmark & \checkmark & & Satisfy Constraints + (above)  \\
       & & \checkmark & Trustworthy Inferences \\
       \rowcolor{\ROWCOLOR} 
       \checkmark & \checkmark & \checkmark & Explainable Errors + (above)
    \end{tabular}   
    \renewcommand{\arraystretch}{1.0}

    \vspace*{5pt}
    \caption{Summary of design features   and corresponding model properties. Design features yield additive properties, as indicated by ``+ (above).'' Proposed implicit L2O models with certificates have intuitive design, memory efficient training, inferences that satisfy optimality/constraint conditions, certificates of trustworthiness, and explainable errors.}
    \label{tab: L2O-comparison}
\end{table}

\subsection{Related Works}
Closely related to our work is deep unrolling, a subset of L2O wherein models consist of a fixed number of iterations of a data-driven optimization algorithm. Deep unrolling has garnered great success and provides intuitive model design. We refer readers to recent surveys \cite{amos2022tutorial,chen2021learning,shlezinger2020model,monga2021algorithm} for further L2O background.
Downsides of unrolling are growing memory requirements with unrolling depth and a lack of guarantees. 

Implicit models circumvent these two shortcomings by defining models using an equation (\eg as in (\ref{eq: L2O-model})) rather than prescribe a fixed number of computations as in deep unrolling. This enables inferences to be computed by iterating until convergence, thereby enabling theoretical guarantees. Memory-efficient training techniques were also developed for this class of models, which have been applied successfully in games, \cite{heaton2021learn} music source separation, \cite{koyama2021music}   language modeling, \cite{bai2019deep}, segmentation, \cite{bai2020multiscale} and inverse problems. \cite{heaton2021feasibility,gilton2021deep}
The recent work \cite{gilton2021deep} most closely aligns with our L2O methodology.

Related XAI works use labels/cards. Model Cards \cite{mitchell2019model}  document intended and appropriate uses of models. Care labels \cite{morik2021care,morik2021yes} are similar, testing properties like expressivity, runtime, and memory usage.
FactSheets \cite{arnold2019factsheets} are modeled after   supplier  declarations of conformity  and aim to identify models' intended use, performance, safety, and security.  
These works provide statistics at the distribution level, complementing our work for trustworthiness of individual inferences.

\section{Explainability via Optimization}

\paragraph{Model Design}
The design of L2O models is naturally decomposed into two steps: optimization formulation and algorithm choice.
The first step is to identify a tentative objective to encode prior knowledge via regularization (\eg sparsity) or constraints (\eg unit simplex for classification). We may also add terms that are entirely data-driven. Informally, this step  identifies a special case of (\ref{eq: L2O-model}) of the form 
\begin{align}
    \begin{split}
        {N}_\Theta(d) \triangleq \argmin_{x} \ \  \mbox{(prior knowledge)}  + \mbox{(data-driven terms)},
    \end{split}
\end{align} 
where the constraints are encoded in the objective using indicator functions, equaling $0$ when constraint is satisfied and $\infty$ otherwise.
The second design step is to choose an algorithm for solving the chosen optimization problem (\eg proximal-gradient or ADMM~ \cite{deng2016global}). 
We use iterative algorithms, and the update formula for each iteration is given by a \textit{model operator} $T_\Theta(x;d)$. Updates are typically composed in terms of gradient and proximal operations.
Some parameters (\eg step sizes) may be included in the weights $\Theta$ to be tuned during training. Given data $d$, computation of the inference ${N}_\Theta(d)$ is completed by generating a sequence $\{x_d^k\}$ via the relation
\begin{equation}
    x_d^{k+1} = T_\Theta(x_d^k; d), \ \ \ \mbox{for all $k\in\bbN$.}
    \label{eq: xk-iteration}
\end{equation}
By design, $\{x_d^k\}$ converges to a solution of (\ref{eq: L2O-model}), and   we set
\begin{equation}
    {N}_\Theta(d)
    = \lim_{k\rightarrow\infty} x_d^k.
\end{equation}  
In our context, each model inference ${N}_\Theta(d)$ is defined to be an optimizer as in (\ref{eq: L2O-model}).  {Hence \textit{properties of inferences can be explained via the optimization model} (\ref{eq: L2O-model}); note this is unlike blackbox models where one has no way of explaining why a particular inference is made.} The iterative algorithm is applied successively until stopping criteria are met (\ie in practice we choose an iterate $K$, possibly dependent on $d$, so that ${N}_\Theta(d) \approx x_d^K$). Because $\{x_d^k\}$ converges, we may adjust stopping criteria to approximate the limit to arbitrary precision, which implies we may provide guarantees on model inferences (\eg satisfying a linear system of equations to a desired precision \cite{gilton2021deep,heaton2021feasibility,heaton2021learn}). The properties of the implicit L2O model (\ref{eq: L2O-model}) are summarized by Table \ref{tab: L2O-comparison}.

\begin{figure*}
    \centering 
    \includegraphics[width=0.30\textwidth]{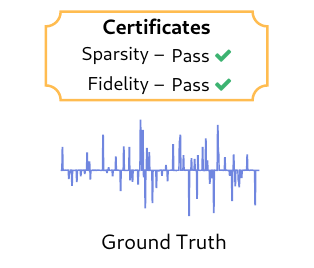} 
    \hspace{-39pt}
    \includegraphics[width=0.30\textwidth]{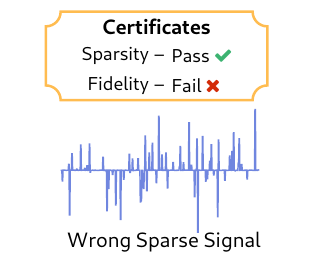}
    \hspace{-39pt}
    \includegraphics[width=0.30\textwidth]{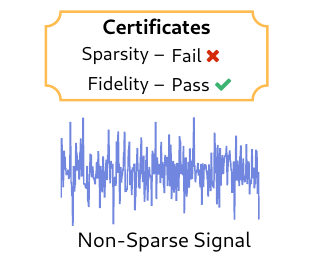}
    \hspace{-39pt}
    \includegraphics[width=0.30\textwidth]{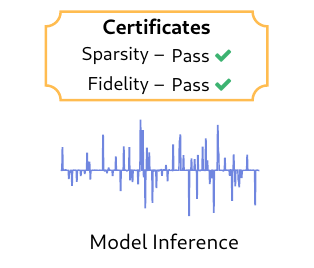}
    
    \vspace{-0.1in}
    
    \caption{Example inferences for test data $d$.  The sparsified version $Kx$ of each inference $x$ is shown (\cf Figure \ref{fig: dictionary-dense-sparse}) along with certificates. 
    Ground truth was taken from test dataset of implicit dictionary experiment. The second from left is  sparse and inconsistent with measurement data. The second from right  complies with measurements but is not sparse. The rightmost   is generated using our proposed model (IDM), which approximates the ground truth well and is trustworthy.     }
    \label{fig: certificate-example}
\end{figure*}

\paragraph{Example of Model Design.} To make the model design procedure concrete, we illustrate this process on a classic problem: sparse recovery from linear measurements. These problems appear in many applications such as radar imaging~ \cite{siddamal2015survey} and speech recognition~ \cite{gemmeke2010compressive}. Here the task is to estimate a signal $x_d^\star$ via access to linear measurements $d$ satisfying $d=Ax_d^\star$ for a known matrix $A$.

\newpage
\textit{Step 1:  {Choose   Model}.} Since true signals are known to be sparse, we include   $\ell_1$ regularization. To comply with measurements, we add a fidelity term. Lastly, to capture hidden features of the data distribution, we also add a data-driven regularization. Putting these together gives the problem 
\begin{equation}
    \min_{x\in\bbR^n} \underbrace{\tau \|x\|_1}_{\mbox{sparsity}} \hspace{-4pt} + \underbrace{\|Ax-d\|_2^2}_{\mbox{fidelity}} + \underbrace{\|W_1Ax\|^2 + \left<x, W_2d\right>}_{\mbox{data-driven regularizer}}\hspace{-1pt},
    \label{eq: lasso}
\end{equation} 
where $\tau > 0$ and $W_1$ and $W_2$ are two tunable matrices. This model encodes a balance of three terms -- sparsity, fidelity, data-driven regularization -- each  quantifiable via (\ref{eq: lasso}).

\textit{Step 2:  {Choose   Algorithm}.}
The proximal-gradient scheme   generates a sequence $\{z^k\}$ converging to a limit which   solves (\ref{eq: lasso}). By simplifying and combining terms,   the proximal-gradient method can be written via the iteration 
\begin{equation}
    z^{k+1} = \eta_{\tau \lambda }\big(z^k - \lambda W (Az^k-d)\big),
    \ \ \ \mbox{for all $k\in\bbN$,}
    \label{eq: proximal-gradient}
\end{equation}
where $\lambda  > 0$ is a step-size, $W$ is a matrix defined in terms
of $W_1$, $W_2$, and $A^\top$, and $\eta_\theta$ is the shrink operator given by 
\begin{equation}
    \eta_{\theta}(x) \triangleq \mbox{sign}(x) \max(|x|-\theta, 0).
\end{equation}
From the update on the right hand side of (\ref{eq: proximal-gradient}), we see the step size $\lambda$ can be ``absorbed'' into the tunable matrix $W$ and the shrink function parameter can be set to $\theta > 0$. That is, this example model has weights $\Theta = (W,\theta, \tau)$ with model operator
\begin{equation}
    T_\Theta(x;d) \triangleq \eta_\theta \big(x - W(Ax-d)\big),
    \label{eq: lasso-T}
\end{equation} 
which resembles the updates of previous L2O works. \cite{liu2019alista,gregor2010learning,chen2018theoretical}
Inferences are computed via a sequence $\{x_d^k\}$ with updates
\begin{equation}
    x_d^{k+1} = T_\Theta(x_d^k; d),
    \ \ \ \mbox{for all $k\in\bbN$.}
    \label{eq: lasso-model-iteration}
\end{equation}
The model inference is  the limit $x_d^\infty$ of this sequence $\{x_d^k\}$.

\paragraph{Convergence} Evaluation of the model ${N}_\Theta(d)$ is well-defined and tractable under a simple assumption. 
By a classic result~ \cite{Krasnoselskii1955_two}, it suffices to ensure, for all $d$, $T_\Theta(\cdot;\ d)$ is \textit{averaged}, \ie there is $\alpha \in (0,1)$ and $Q$ such that 
$T_\Theta(x;d) = (1-\alpha)x + \alpha Q(x;d)$, where $Q$ is $1$-Lipschitz in $x$.  
When this property holds, the sequence $\{x_d^k\}$ in (\ref{eq: xk-iteration}) converges to a solution $x^\star_d$. This may appear to be a strong assumption; however, common operations in  convex optimization algorithms (\eg proximals and gradient descent updates) are averaged. For entirely data-driven portions of $T_\Theta$,   several activation functions are  1-Lipschitz \cite{combettes2020lipschitz,gao2017properties} (\eg ReLU and softmax), and libraries like PyTorch \cite{paszke2019pytorch} include functionality to force affine mappings to be 1-Lipschitz (\eg spectral normalization). 
Furthermore, by making $T_\Theta(\cdot;d)$ a contraction, a unique fixed point is obtained.
We emphasize, even without forcing $T_\Theta$ to be averaged,   $\{x^k\}$ is often   observed to converge in practice \cite{heaton2021feasibility,bai2019deep,gilton2021deep} upon tuning   weights $\Theta$. 

\subsection{Trustworthiness Certificates} 
Explainable models justify whether each inference is trustworthy. 
We propose providing justification in the form of certificates, which verify various properties of the inference are consistent with those of the model inferences on training data and/or prior knowledge. Each certificate is a tuple of the form $(\mbox{name}, \mbox{label})$ with a property name and a corresponding label which has one of three values: pass, warning, or fail (see Figure~\ref{fig: certificate-example}).
Each certificate label is generated by two steps. The first is to apply a function that maps inferences (or intermediate states) to a \textit{nonnegative scalar value $\alpha$} quantifying a property of interest. The second step is to map this scalar to a label. Labels are generated via the flow:
\begin{equation}
    (\mbox{Inference})
    \rightarrow   
    (\mbox{Property Value})
    \rightarrow  
    (\mbox{Certificate Label}).
\end{equation} 

\paragraph{Property Value Functions} 
Several quantities may be used to generate certificates. 
In the model design example above, a sparsity property can be quantified by counting the number of nonzero entries in a signal, and a fidelity property can use the relative error $\|Ax-d\|/\|d\|$ (see Figure \ref{fig: certificate-example}).
To be most effective, property values are chosen to coincide with the optimization problem used to design the L2O model, \ie to quantify structure of prior and data-driven knowledge. 
This enables each certificate to clearly validate a portion of the model (see Figure \ref{fig: certificate-diagram}). Since various concepts are useful for different types of modeling, we provide a brief (and non-comprehensive) list of concepts and possible corresponding property values in Table \ref{tab: certificate-diagonistics}. 

\begin{table}[t]
    \centering 
    \renewcommand{\arraystretch}{1.25}
    \begin{tabular}{c|c|c}
        Concept & Quantity & Formula \\\hline 
        
        Sparsity & Nonzeros  & $\|x\|_0$ \\

        \rowcolor{\ROWCOLOR}
        $\approx$ Sparsity & $\ell_1$ norm & $\|x\|_1$\\

        Measurements & Relative Error & $\|Ax-d\|/\|d\|$\\ 

        \rowcolor{\ROWCOLOR}
        Soft Constraint & Distance to Set $\sC$ & $d_{\sC}(x)$ \\

        Smooth Images & Total Variation & $\|\nabla x\|_1$\\

        \rowcolor{\ROWCOLOR}
        Convergence & Iterate Residual & $\|x^{k}-x^{k-1}\|$\\
        
        Classifier   & Probability short & 
         \multirow{2}*{$1- \max_i x_i$} \\[-0.5pt] 
        Confidence & of one-hot label &      \\ 
        
        \rowcolor{\ROWCOLOR}
        Regularization & Proximal Residual & $\|x-\prox{f_\Omega}(x)\|$ 
    \end{tabular} 
    \renewcommand{\arraystretch}{1.0} 
    \vspace*{5pt}
    
    \caption{Each certificate is tied to a high-level concept, and its corresponding property value is quantified using a formula. For classifier confidence, we assume $x$ is in the unit simplex. The proximal is a data-driven update for $f_\Omega$ with weights $\Omega$.}
    \label{tab: certificate-diagonistics}
\end{table}

One property concept deserves particular attention: data-driven regularization.  This regularization is important for discriminating between inference features that are qualitatively intuitive but  difficult to quantify by hand. 
Rather than approximate a function, implicit L2O models directly approximate gradients/proximals.  
These provide a way to measure regularization indirectly via gradient norms/residual norms of proximals. 
Moreover, these norms (\eg see last row of Table 2) are easy to compute and equal zero  only at  local minima of   regularizers. 
To our knowledge, this is the first work to \textit{quantify} trustworthiness using the quality of inferences with respect to data-driven regularization.

\paragraph{Certificate Labels}
Typical certificate labels should follow a trend where inferences often obtain a pass label to indicate trustworthiness while warnings occur occasionally and failures are obtained in extreme situations.  
Let the samples of model inference property values $\alpha \in [0,\infty)$ come from distribution $\bbP_\sA$.
We pick property value functions for which small $\alpha$ values are desirable and the distribution tail consists of larger $\alpha$.
Intuitively, smaller property values of $\alpha$ resemble property values of inferences from training and/or test data. 
Thus, labels are assigned according to the probability of observing a value less than or equal to $\alpha$, \ie  
we evaluate the cumulative distribution function (CDF) defined for probability measure $\bbP_{\sA}$ by 
\begin{equation}
    \mathrm{CDF}(\alpha) = \int_0^\alpha \ \mathrm{d} \bbP_\sA,
\end{equation}
\begin{figure}
    \centering
    \includegraphics{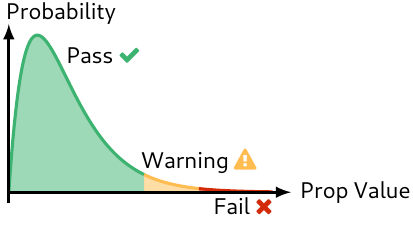}
    \caption{Probability distribution for values and labels of a particular model property. The majority of samples drawn from this distribution are set to pass while the outliers in the tail fail.}
    \label{fig: cert-prob}
\end{figure}
Labels are chosen according to the task at hand. 
Let $p_{\mathrm{p}}$, $p_{\mathrm{w}}$, and $p_{\mathrm{f}} = 1- p_{\mathrm{p}}-p_{\mathrm{w}}$ be the   probabilities for   pass, warning, and fail labels, respectively.
Labels are made for $\alpha$ via
\begin{equation}
    \mbox{Label}(\alpha) =
    \begin{cases}
    \begin{array}{cl}
        \mbox{pass}     &  \mbox{if $\mathrm{CDF}(\alpha) < p_{\mathrm{p}}$ } \\
        \mbox{warning}  & \mbox{if $\mbox{CDF}(\alpha) \in [p_{\mathrm{p}}, 1- p_{\mathrm{f}})$}\\
        \mbox{fail}    &  \mbox{otherwise}.    
    \end{array}
    \end{cases}
    \label{eq: label_assignment}
\end{equation} 

The remaining task is to estimate the CDF value for a given $\alpha$.
Recall we assume   access is given to property values $\{\alpha_i\}_{i=1}^N$ from ground truths or inferences on training data, where $N$ is the number of data points.
To this end, given an $\alpha$ value, we estimate its CDF value via the empirical CDF: 
\begin{align}
    \mathrm{CDF}(\alpha) 
    \approx 
    \dfrac{|\{ \alpha_i : \alpha_i \leq \alpha,\ 1\leq i \leq N\}|}{N} 
    = \dfrac{\# \text{ of } \text{$\alpha_i$'s} \leq \alpha}{N},
    \label{eq: cdf}
\end{align}
where $|\cdot|$ denotes set cardinality. 
For large $N$, (\ref{eq: cdf})   well approximates the continuous CDF.

\paragraph{Certificate Implementation} 
As noted in the introduction, trustworthiness certificates are evidence an inference satisfies post-conditions (\ie passes various tests).   Thus, they are to be used in code in the same manner as standard software engineering practice. Consider the snippet of code in Figure \ref{fig: python-sample-code}. As usual, an inference is generated by calling the model. However, alongside the inference \verb|x|, certificates \verb|certs| are returned that label whether the inference \verb|x| passes tests that identify consistency with training data and prior knowledge.

\section{Experiments}
Each numerical experiment shows an application of novel implicit L2O models, which were designed directly from prior knowledge. Associated certificates of trustworthiness are used to emphasize the explainability of each model and illustrate use-cases of certificates. Experiments were coded using Python with the PyTorch library  \cite{paszke2019pytorch}, the Adam optimizer \cite{kingma2014adam}, and, for ease of re-use, were run via Google Colab. 
We emphasize these experiments are for illustration of intuitive and novel model design and trustworthiness and are not    benchmarked against state-of-the-art models.
The datasets generated and/or analysed during the current study are available in the following repository: \href{https://github.com/typal-research/xai-l2o}{github.com/typal-research/xai-l2o}.
All methods were performed in accordance with the relevant guidelines and regulations.

{\subsection{Algorithms}
To illustrate   evaluation of  L2O model used herein, we begin with an example  L2O model and algorithm. Specifically,   models used for the first two experiments take the form
\begin{equation}
    \min_{x\in\bbR^n} f(Kx) + h(x)
    \ \ \mbox{s.t.}\ \ 
    \|Mx-d\|\leq \delta, 
\end{equation}
where $K$ and $M$ are linear operators, $\delta \geq 0$ is a noise tolerance, and $f$ and $h$ are proximable\footnote{A function is proximable if it admits a ``nice'' closed-form proximal formula, where $\prox{f}(x)\triangleq \underset{z}{\argmin} f(z) + \|z-x\|^2$.} functions. 
Introducing auxiliary variables $w$ and $p$ and dual variable $\nu=(\nu_1,\nu_2)$, linearized ADMM~ \cite{ryu2022large} (L-ADMM) can be used to iteratively update the tuple $(p,w,\nu, x)$ of variables via
{\small 
\begin{subequations}
\begin{align}
    p^{k+1} & = \prox{\lambda f} \left(p^k + \lambda(\nu_1^k + \alpha(Kx^k - p^k)) \right)
    \\
    w^{k+1} & = \mbox{proj}_{B(d,\delta)}\left(w^k + \lambda(\nu_2^k + \alpha(Mx^k - w^k)) \right)
    \\
    \nu_1^{k+1} &= \nu_1^k + \alpha(Kx^k - p^{k+1})
    \\
    \nu_2^{k+1} &= \nu_2^k + \alpha(Mx^k - w^{k+1})
    \\
    r^{k} & =  K^\top \left( 2\nu_1^{k+1} - \nu_1^k  \right) + M^\top \left( 2\nu_2^{k+1} - \nu_2^k \right)\\
    x^{k+1} & =\prox{\beta h}\left(x^k - \beta r^k\right),
\end{align}\label{eq: L-ADMM-main-text-form}\end{subequations}}where $\mbox{proj}_{B(d,\delta)}$ is the Euclidean projection onto the Euclidean ball of radius $\delta$ centered at $d$,  $\prox{f}$ is the proximal operator for a function $f$, and the scalars $\alpha,\beta,\lambda > 0$ are appropriate step sizes. 
Further details, definitions, and explanations are available in the appendices. Note the updates are ordered so that $x^{k+1}$ is the final step to make it easy to backprop through the final $x^k$ update.}

\subsection{Implicit Model Training}
Standard backpropagation cannot be used for implicit models  as it requires memory capacities beyond   existing computing devices. Indeed,  storing gradient data for each iteration in the forward propagation (see (\ref{eq: xk-iteration})) scales the memory  during training linearly with respect to the number of iterations.
Since the limit $x^\infty$ solves a fixed point equation, implicit models can be trained by differentiating implicitly through the fixed point to obtain a gradient. This implicit differentiation requires further computations and coding. Instead of using gradients, we utilize Jacobian-Free Backpropagation (JFB) \cite{fung2021jfb} to train models. JFB further simplifies training by only backpropagating through the final iteration, which was proven to yield  preconditioned gradients. JFB trains using fixed memory  (with respect to the $K$ steps used to estimate ${N}_\Theta(d)$) and avoids numerical issues   arising from  computing exact gradients, \cite{bai2021stabilizing} making JFB and its variations \cite{geng2021training,huang2021textrm} apt for training implicit models.

\subsection{Implicit Dictionary Learning} \label{subsec: exp-dict}

\paragraph{Setup} 
In practice, high dimensional signals often approximately admit low dimensional representations~ \cite{osher2017low,zhang2015survey,carlsson2008local,lee2003nonlinear,peyre2008image,peyre2009manifold}.
For illustration, we   consider a linear inverse problem where true data admit sparse representations. 
Here each signal $x_d^\star \in\bbR^{250}$ admits a   representation  $s_d^\star\in\bbR^{50}$ via a  transformation $M$ (\ie $x_d^\star = Ms_d^\star$).  A matrix  $A\in\bbR^{100\times250}$ is applied to each signal  $x_d^\star$  to provide   linear measurements $d = Ax_d^\star$. Our task is to recover $x^\star_d$ given knowledge of $A$ and $d$ \textit{without} the   matrix $M$.   Since the linear system is quite under-determined, schemes solely minimizing measurement error  (\eg  least squares approaches)  fail to recover true signals; additional knowledge is essential. 

\begin{figure}[t] 
    \centering 
    \includegraphics[width=3.0in]{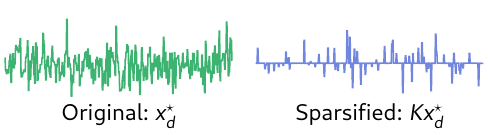}   
    \vspace{-0.1in}
    \caption{Training IDM yields sparse representation of inferences. Diagram shows a sample true data $x$ (left) from test dataset   and its sparsified representation $Kx$  (right). }
    \label{fig: dictionary-dense-sparse}
\end{figure}

  \begin{figure*}
    \def\spycolor{red!75!black}
    \def\W{1.5pt}
      \centering
      \small
        \setlength{\tabcolsep}{0.1pt}
        \begin{tabular}{ccccc}

        \begin{tikzpicture} [spy using outlines={rectangle, magnification=3, size=1cm, connect spies}, rounded corners]
                \node[anchor=south west,inner sep=0] (image) at (0,0) {\adjincludegraphics[width=0.185\textwidth, rotate=180]{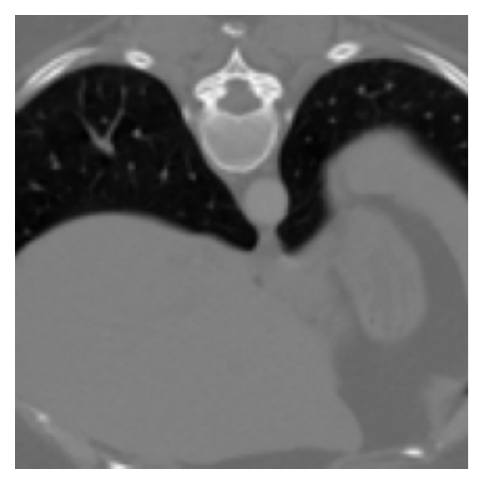}}; 
                \begin{scope}
                \spy[\spycolor,size=0.175\textwidth, every spy on node/.append style={line width = \W}] on (1.65,0.9) in node at (1.65, -1.55);  
                \end{scope}
        \end{tikzpicture}
        &
        \begin{tikzpicture} [spy using outlines={rectangle, magnification=3, size=1cm, connect spies}, rounded corners]
                \node[anchor=south west,inner sep=0] (image) at (0,0) {\adjincludegraphics[width=0.185\textwidth, rotate=180]{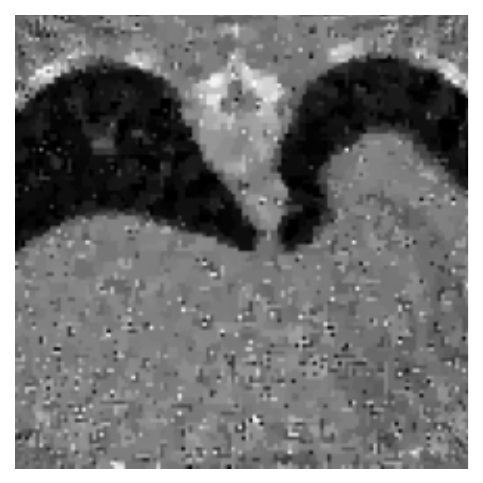}};
                \begin{scope}
                \spy[\spycolor,size=0.175\textwidth, every spy on node/.append style={line width = \W}] on (1.65,0.9)  in node at (1.65, -1.55);  
                \end{scope}
        \end{tikzpicture}   
        &
        \begin{tikzpicture} [spy using outlines={rectangle, magnification=3, size=1cm, connect spies}, rounded corners]
                \node[anchor=south west,inner sep=0] (image) at (0,0) {\adjincludegraphics[width=0.185\textwidth, rotate=180]{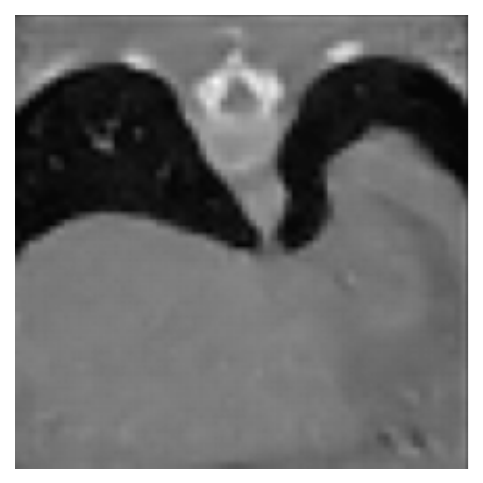}};
                \begin{scope}
                \spy[\spycolor,size=0.175\textwidth, every spy on node/.append style={line width = \W}] on (1.65,0.9) in node at (1.65, -1.55);  
                \end{scope}
        \end{tikzpicture}       
        &
        \begin{tikzpicture} [spy using outlines={rectangle, magnification=3, size=1cm, connect spies}, rounded corners]
                \node[anchor=south west,inner sep=0] (image) at (0,0) {\adjincludegraphics[width=0.185\textwidth, rotate=180]{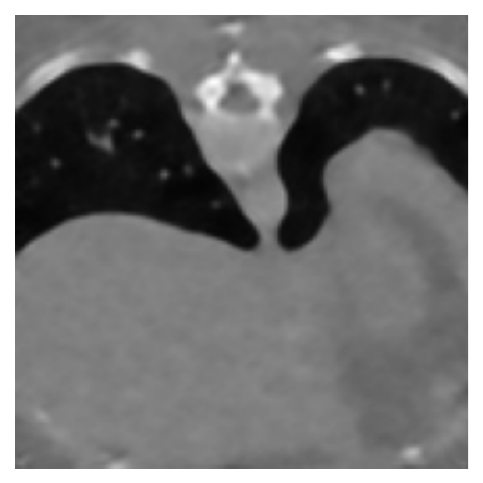}};
                \begin{scope}
                \spy[\spycolor,size=0.175\textwidth, every spy on node/.append style={line width = \W}]on (1.65,0.9)  in node at (1.65, -1.55);  
                \end{scope}
        \end{tikzpicture}       
        &
        \begin{tikzpicture} [spy using outlines={rectangle, magnification=3, size=1cm, connect spies}, rounded corners]
                \node[anchor=south west,inner sep=0] (image) at (0,0) {\adjincludegraphics[width=0.185\textwidth, rotate=180]{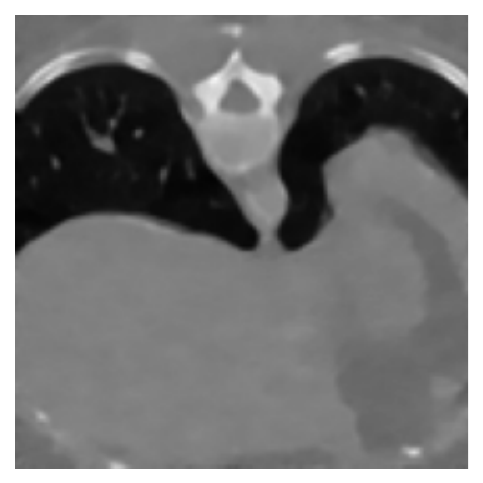}};
                \begin{scope}
                \spy[\spycolor,size=0.175\textwidth, every spy on node/.append style={line width = \W}] on (1.65,0.9)  in node at (1.65, -1.55);  
                \end{scope}
        \end{tikzpicture}   
        \\
        {\bf Ground truth} & {\bf TV Min} & {\bf U-Net} & {\bf F-FPN} & {\bf Implicit L2O}
        \\      
         
        SSIM: 1.000 & SSIM: 0.379 & SSIM: 0.806 & SSIM: 0.853 & SSIM: 0.891 \\ 
        PSNR: $\infty$ & PSNR: 21.17 & PSNR: 27.83 & PSNR: 30.32 & PSNR: 31.87 \\  
        & 
        {Fail \textcolor{RED}{\faTimes} -- Data Regularization}
        & 
        {Fail \textcolor{RED}{\faTimes} -- Fidelity}
        & 
        {Fail \textcolor{RED}{\faTimes} -- Box Constraint}
        & 
        {Trustworthy \textcolor{GREEN}{\faCheck}} \\   
        \end{tabular}
        \caption{Reconstructions on test data computed via U-Net, \cite{jin2017deep} TV minimization, F-FPNs, \cite{heaton2021feasibility} and Implicit L2O (left to right). Bottom row shows expansion of region indicated by red box. Pixel values outside $[0,1]$ are flagged.  Fidelity is flagged when images do not comply with measurements, and regularization is flagged when texture features of images are sufficiently inconsistent with true data (\eg grainy images). Labels are provided beneath each image (\nb fail is assigned   to images that are worse \textcolor{black}{than 95\% of L2O inferences on training data}). Shown comparison methods fail  while the Implicit L2O image passes all   tests.}
        \label{fig: LoDoPab-Examples}
    \end{figure*}

\paragraph{Model Design}  
All convex regularization approaches are known lead to biased estimators whose expectation does not equal the true signal. \cite{fan2001variable}
However, the seminal work \cite{candes2006robust} of Candes and Tao shows $\ell_1$ minimization (rather than additive regularization) enables exact recovery under suitable assumptions.
Thus, we minimize a sparsified signal subject to linear constraints via the implicit dictionary model (IDM)
\begin{equation}
    {N}_\Theta(d) \triangleq   \argmin_{x\in\bbR^{250}}    \|Kx\|_1 \ \ \mbox{s.t.} \ \ Ax = d.
    \label{eq: model-sparse-recovery}
\end{equation}  
The square matrix $K$ is used to leverage the fact $x$ has a low-dimensional representation by transforming $x$ into a sparse vector.  Linearized ADMM  \cite{ryu2022large} (L-ADMM) is used to create a sequence $\{x_d^k\}$ as in (\ref{eq: xk-iteration}). The model ${N}_\Theta$ has weights $\Theta = K$.
If it exists, the matrix $K^{-1}$ is known as a dictionary and $K {N}_\Theta(d)$ is the corresponding sparse code; hence the name IDM for (\ref{eq: model-sparse-recovery}).
To this end, we emphasize $K$ is learned during training and is \textit{different} from $M$, but these matrices are related since we aim for the product
$
    Kx_d^\star
    = KMs_d^\star
$
to be sparse. 
Note we use  L-ADMM  to  \textit{provably} solve  (\ref{eq: model-sparse-recovery}), and ${N}_\Theta$ is easy to train.
More details can be found in the appendix.

\paragraph{Discussion}
IDM combines intuition from dictionary learning with a  reconstruction algorithm. Two properties are used to identify trustworthy inferences: sparsity and measurement compliance (\ie fidelity). Sparsity and fidelity are  quantified via the $\ell_1$ norm of the sparsified inference (\ie $K{N}_\Theta(d)$) and  relative measurement error. Figure \ref{fig: dictionary-dense-sparse} shows the training the model yields a sparsifying transformation $K$.  
 Figure \ref{fig: certificate-example} shows the proposed certificates   identify   ``bad''   inferences that might, at first glance, appear to be ``good'' due to their compatibility with constraints. 
Lastly, observe the utility of learning $K$, rather than  approximating $M$, is  $K$ makes it is easy to check if an inference admits a sparse representation. Using $M$ to check for sparsity is nontrivial.
 
\subsection{CT Image Reconstruction}\label{subsec: exp-ct} 
    \paragraph{Setup} Comparisons are provided for low-dose CT examples derived from the Low-Dose Parallel Beam dataset (LoDoPab) dataset, \cite{leuschner2019lodopab} which has publically available phantoms  derived from actual human chest CT scans.   CT measurements are simulated with a parallel beam geometry and a sparse-angle setup of only $30$ angles and $183$ projection beams, giving   5,490 equations and 16,384 unknowns. We add $1.5\%$ Gaussian noise to \emph{each individual beam measurement}. Images have resolution $128 \times 128$.
     To make errors  easier to contrast between methods,  the linear systems here are under-determined and have more noise  than those in some similar works. Image quality is determined using the Peak  Signal-To-Noise Ratio (PSNR) and structural similarity index measure
    (SSIM).  The training loss was mean squared error.  Training/test datasets have 20,000/2,000 samples.

    \paragraph{Model Design} 
    The model for the CT experiment extends the IDM. In practice, it has been helpful to utilize a sparsifying transform. \cite{jiang2018super,xu2012low} We accomplish this via a linear  operator $K$, which is applied and then this product is fed into a data-driven regularizer $f_\Omega$ with parameters $\Omega$. We additionally ensure compliance with measurements from the Radon transform matrix $A$, up to a tolerance $\delta$. In our setting, all pixel values are also known to be in the interval $[0,1]$. Combining our prior knowledge yields the implicit L2O model 
    \begin{equation}
      {N}_\Theta(d)
      \triangleq \argmin_{x\in[0,1]^n} f_\Omega(Kx) \ \ \mbox{s.t.}\ \ \|Ax-d\| \leq \delta.
      \label{eq: model-ct}
    \end{equation}  
    Here ${N}_\Theta$ has weights $\Theta = (\Omega, K,  \alpha, \beta, \lambda)$ with $\alpha$, $\beta$ and $\lambda$ step-sizes in L-ADMM. 
    More details can be found in the appendix.

    \begin{table*}[t] 
        \centering 
        \renewcommand{\arraystretch}{1.25}
        \begin{tabular}{c|c|c|c|c|c|c}
        {Method} & Avg. PSNR   & Avg. SSIM & Box Constraint Fail & Fidelity Fail & Data Reg. Fail & \# Params 
        \\\hline
        \rowcolor{\ROWCOLOR}
        U-Net & 27.32 dB & 0.761 & 5.75 \% & 96.95\% &  3.20\% & 533,593
        \\          
        TV Min & 28.52 dB & 0.765 &  0.00 \% & 0.00\% &  25.40\% & 4
        \\          
        \rowcolor{\ROWCOLOR}
        F-FPN\textsuperscript{$\dagger$} & 30.46 dB & 0.832 & 47.15\% &  0.40\% & 5.05\% & 96,307
        \\   
        Implicit L2O   & 31.73 dB & 0.858 & 0.00\% & 0.00\% &  5.70\% & 59,697
        \end{tabular}
        \renewcommand{\arraystretch}{1.0}
        \caption{Average PSNR/SSIM for CT reconstructions on the 2,000 image LoDoPab testing dataset. $\dagger$ Reported from original work. \cite{heaton2021feasibility}  U-Net was trained with filtered backprojection as in prior work. \cite{jin2017deep}
        Three properties are used to check trustworthiness: box constraints, compliance with measurements (\ie  fidelity), and data-driven regularization   (via the proximal residual   in Table~\ref{tab: certificate-diagonistics}). Failed sample percentages  are numerically estimated via~\eqref{eq: cdf}. Sample property values ``fail'' if they perform worse than $95\%$ of the inferences on the training data, \ie  its CDF value exceeds 0.95. Implicit L2O yields the most passes on test data.}
        \label{tab: LoDoPab_results}
        \vspace*{-0.2in}
    \end{table*}    
   
   \paragraph{Discussion}
   Comparisons of our method (Implicit L2O) with U-Net, \cite{jin2017deep} F-FPNs, \cite{heaton2021feasibility} and total variation (TV) Minimization are given in Figure \ref{fig: LoDoPab-Examples} and Table \ref{tab: LoDoPab_results}. 
   Table~\ref{tab: LoDoPab_results} shows the average PSNR and SSIM reconstructions. Our model obtains the highest average PSNR and SSIM values on the test data while using 11\% and 62\% as many weights as U-Net and F-FFPN, indicating greater efficiency of the implicit L2O framework.
   Moreover, the L2O model is designed with three features: compliance with measurements (\ie fidelity), valid pixel values, and data-driven regularization.
   Table~\ref{tab: LoDoPab_results} also shows the percentage of ``fail'' labels for these property values.
   Here,   an inference fails if its property value is larger than $95\%$ of the property values from the training/true data, \ie  we choose $p_p = 0.95$, $p_w = 0$, and $p_f = 0.05$ in  \eqref{eq: label_assignment}.
   For the fidelity,  our model never fails (due to incorporating the constraint into the network design). 
   Our network fails \textcolor{black}{$5.7\%$} of the time for the data-driven regularization property. 
   Overall, the L2O model generates the most trustworthy inferences.
   This is intuitive as this  model  outperforms the others and was specifically designed to embed all of our knowledge, unlike the others.   
    To provide better intuition of the certificates, we also show the certificate labels for an image from the test dataset in  Figure \ref{fig: LoDoPab-Examples}.
    The only image to pass all  provided tests is the proposed implicit L2O model.  
   This knowledge can help identify trustworthy inferences. 
   Interestingly,  the data-driven regularization enabled   certificates to detect and flag ``bad'' TV Minimization features (\eg  visible staircasing effects \cite{ring2000structural,chan2000high}), which shows novelty of certificates as these features are intuitive, yet  prior methods to quantify this were, to our knowledge, unknown.


\subsection{Optimal Cryptoasset Trading}
\paragraph{Setup} 
Ethereum is a blockchain technology    anyone can use to deploy permanent and immutable decentralized applications. 
This technology enables creation of decentralized finance (DeFi) primitives, which can give censorship-resistant participation in   digital markets and expand the use of   stable assets \cite{kamvar2019celo,makerdao2020makerdao}  and exchanges \cite{zhang2018formal,warren20170x,hertzog2017bancor}   beyond the realm of traditional finance. Popularity of  cryptoasset trading (\eg GRT and Ether) is exploding   with the DeFi movement. \cite{werner2021sok,schar2021decentralized}

Decentralized exchanges (DEXs) are a popular entity for exchanging cryptoassets (subject to a small transaction fee), where trades  are conducted without the need for a trusted intermediary to facilitate the exchange. 
Popular examples of DEXs are constant function market makers (CFMMs), \cite{angeris2020improved} which use mathematical formulas to govern trades. To ensure   CFMMs maintain  sufficient net assets,
trades within   CFMMs maintain constant total reserves (as defined by a function $\phi$).   
A transaction in a CFMM tendering $x$ assets in return for $y$ assets   with reserves assets $r$ is accepted provided 
\begin{equation}
    \phi( r + \gamma x - y) \geq \phi(r),
    \label{eq: CFMM-transaction-rule}
\end{equation}
with $\gamma\in(0,1]$   a trade fee parameter.
Here $r,x,y\in\bbR^n$ with each vector nonnegative and  $i$-th entry giving an amount for the $i$-th cryptoasset type   (\eg Ether, GRT).
Typical choices \cite{angeris2021constant} of $\phi$ are weighted sums and products, \ie 
\begin{equation}
    \phi(r) = \sum_{i=1}^n w_i r_i
    \ \ \mbox{and} \ \ 
    \phi(r) = \prod_{i=1}^n r_i^{w_i}    
    .
\end{equation}
where $w\in\bbR^n$ has positive entries.  

This experiment aims to maximize arbitrage. 
Arbitrage is the  simultaneous purchase and sale of equivalent assets in multiple markets to exploit price discrepancies between the markets.
This can be a lucrative endeavor with cryptoassets. \cite{makarov2020trading}
For a given snapshot in time, our arbitrage goal is to identify a collection of trades that maximize the cryptoassets  obtainable by   trading between different exchanges, \ie solve the (informal) optimization  problem  
\begin{equation}
    \max_{\small \mathrm{trade}} \mbox{Assets}(\mbox{trade})
    \ \ \mbox{s.t.} \ \ 
    \mbox{trade}\in \{\mbox{valid trades}\}.
    \label{eq: model-crypto-informal}
\end{equation}
The set of valid trades is all trades satisfying the transaction rules for CFMMs given by (\ref{eq: CFMM-transaction-rule})   with nonnegative values for tokens tendered and received (\ie $x,y\geq 0$). 
Prior works \cite{angeris2021optimal,angeris2021constant} deal with an idealistic noiseless setting while recognizing executing trades is not without risk (\eg  noisy information, front running, \cite{daian2019flash}  and trade delays). 
To show  implications of trade risk, we incorporate noise in our trade simulations by adding noise $\varepsilon \in \bbR^n$ to CFMM asset observations, which yields noisy observed data $d = (1 +  \varepsilon) \odot r$. 
Also,   we consider trades with CFMMs where \textit{several} assets can be traded simultaneously rather than restricting to pairwise  swaps.

\begin{figure}
    \centering
    \includegraphics{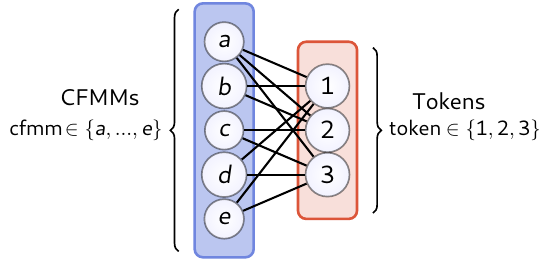}
    \caption{Network with 5 CFMMs and 3 tokens; structure   replicates an experiment in recent work. \cite{angeris2021optimal} Black lines show  available tokens for trade in each CFMM.}
    \label{fig: cfmm-token-diagram}
\end{figure}

 \begin{table*}[t] 
    \small 
    \centering 
    \renewcommand{\arraystretch}{1.25}
    \begin{tabular}{c|c|c|c|c|c|c} 
    {Method} & Predicted Utility 
    &  Executed Utility
    &   Trade Execution 
    & Risk Fail
    & Profitable Fail
    &  \# Params 
    \\\hline
    \rowcolor{\ROWCOLOR} 
    Analytic  & 11.446  & 0.00 & 0.00\% &  100.00\% & 0\% & 0
    \\    
    Implicit L2O   & 0.665 &  0.6785 & 88.20\% &  3.6\% &  11.80\% & 126
    \end{tabular}
    \renewcommand{\arraystretch}{1.0}
    \caption{Averaged results on test data for trades in CFMM network. The analytic method always predicts a profitable trade, but fails to satisfy the constraints (due to noise). 
    This failure is predicted by the certificates ``risk'' certificate and reflected by the 0\% trade execution.  Alternatively, the L2O scheme makes conservative predictions regarding constraints, which limits profitability.
    However, using these certificates, executed L2O trades are always profitable and satisfy constraints. 
    }
    \label{tab: crypto_results}
    \vspace*{-0.2in}
\end{table*}

\paragraph{Model Design}
The aim is to create a model that infers a trade $(x,y)$   maximizing   utility. For a nonnegative vector $p\in\bbR^n$ of reference price valuations, this utility $U$ is the net change in asset values provided by the trade, \ie  \\ 
\begin{equation}
    U(x,y) \triangleq  \underbrace{
    \sum_{j=1}^m  \left<A^jp, A^j( y^j - x^j)\right>}_{\mbox{{ net asset value change}}},
\end{equation}
where $A^j$ is a matrix   mapping global coordinates of asset vector to the coordinates of the $j$-th CFMM (see appendix for details).
For noisy data $d$, trade predictions can include  a ``cost of risk.'' This is quantified by regularizing the trade utility, \ie   introducing a penalty term.
For  matrices $W^j$, we model risk by  a  simple quadratic penalty via  \\[-10pt]
\begin{equation}
    U_\Theta(x,y) \triangleq U(x,y)  - 
    \underbrace{\dfrac{1}{2}\hspace{-1.5pt}\cdot\hspace*{-1.0pt} \sum_{j=1}^m \| A^j W^j(x-y)\|^2.
    \hspace*{-2pt}}_{\mbox{risk model}} 
\end{equation} 
The implicit L2O model infers optimal trades via $U_\Theta$, \ie 
\begin{equation} 
    {N}_\Theta(d) \triangleq (x_d,y_d) = \argmax_{(x,y)\in\sC_\Theta(d)} U_\Theta(x,y),
    \label{eq: model-crypto}
\end{equation}
where  $\sC_\Theta(d)$ encodes   constraints for valid transactions. The essence of ${N}_\Theta$ is to output  solutions to (\ref{eq: model-crypto-informal}) that account for transaction risks.  A formulation of Davis-Yin operator splitting \cite{davis2017three} is used for model evaluation.  
Further details of the optimization scheme can be found in the appendix.

 \begin{figure}
    \centering 

    \includegraphics[]{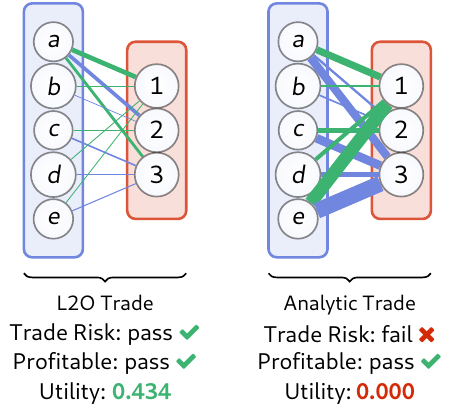}
    
    \caption{Example of proposed L2O (left) and analytic (right) trades with noisy data $d$. Blue and green lines show proposed cryptoassets $x$ and $y$ to tender and receive, respectively (widths show magnitude). The analytic trade  is unable to account for trade risks, causing it to propose large trades that are \textit{not} executed (giving executed utility of zero). 
    This can be anticipated by the failed trade risk certificate. On the other hand, the L2O scheme is profitable (utility is 0.434) and is executed (consistent with the pass trade risk label). 
    }
    \label{fig: cfmm-trade-diagram-01}
\end{figure}
 \newpage

\paragraph{Discussion}
The L2O model contains three core features: profit, risk, and trade constraints. 
The model is designed to output trades that satisfy provided constraints, but note these are \textit{noisy} and thus cannot be used to a priori determine whether a trade will be executed. For this reason,  fail flags identify conditions to warn a trader when a trade should be aborted (due to an ``invalid trade''). 
This can avoid wasting transaction fees (\ie gas costs).
Figure \ref{fig: cfmm-trade-diagram-01} shows an example of two trades, where we note the analytic method proposes a large trade that is \textit{not} executed since it violates the trade constraints (due to noisy observations). The L2O method proposes a small trade that yielded arbitrage profits (\ie $U > 0$) and has pass certificates. 
Comparisons are provided in Table \ref{tab: crypto_results} between the analytic and L2O models. Although the  analytic method has ``ideal'' structure, it performs much worse than the L2O scheme. In particular, \textit{no trades} are executable by the analytic scheme since the present noise always makes the proposed transactions fail to satisfy the actual CFMM constraints. Consistent {with} this, every proposed trade by the analytic trade is flagged as risky in Table \ref{tab: crypto_results}. The  noise is  on the order of {0.2\% Gaussian noise} of the asset totals.

\section{Conclusions} 
Explainable ML models can be concretely developed by fusing certificates with the L2O methodology. 
The implicit L2O methodology enables  prior and data-driven knowledge to be  directly embedded into   models, thereby providing clear and intuitive design. This approach is theoretically sound and compatible with state-of-the-art ML tools. 
The L2O model also enables construction of our  certificate framework  with easy-to-read labels,   certifying if each inference is trustworthy.  
In particular, our certificates provide a principled scheme for the detection  of inferences with ``bad'' features via data-driven regularization. 
Thanks to this optimization-based model design (where inferences can be defined by fixed point conditions), failed certificates can be used to discard untrustworthy inferences and may help debugging the architecture. This reveals the interwoven nature of pairing  implicit L2O with certificates. 
Our experiments  illustrate these ideas in three different settings,  presenting novel model designs and interpretable results.
Future work will study extensions to physics-based applications where PDE-based physics can be integrated into the model  \cite{raissi2019physics,ruthotto2020machine,lin2021alternating}.

\section*{Acknowledgments}
The authors thank Wotao Yin, Stanley Osher, Daniel McKenzie, Qiuwei Li, and Luis Tenorio for many fruitful discussions.  Howard Heaton and Samy Wu Fung were supported by AFOSR MURI FA9550-18-1-0502 and ONR grants: N00014-18-1-2527, N00014-20-1-2093, and N00014-20-1-2787. Samy Wu Fung was partially funded by the National Science Foundation award DMS-2309810.

\printbibliography

\appendix

\section{{Certificate Code Snippet}}

{The figure below shows how code for a \texttt{model} function can be tweaked to include certificates per standard software engineering practice.}

\begin{figure}[h]
\begin{tcolorbox}[colback=white!96!black,colframe=black!20!white,title=Inference + Certificates $\rightarrow$ Trustworthy Inference, coltitle=black, left=1mm]
\vspace*{-2pt}
\small 
\begin{lstlisting}[language=Python, numbers=none, gobble=0,] 
def TrustworthyInference(d):
    x, certs = model(d)
    if 'warning' in certs:
        warnings.warn('Warning Msg')    
    if 'fail' in certs:
        raise Exception('Error Msg')
    return x
\end{lstlisting} 
\end{tcolorbox} 
    \caption{{Example Python code to use certificates as post-conditions. Actual code should use specific warning/exception messages for flagged entries in \texttt{certs}.}  }
    \label{fig: python-sample-code}
\end{figure}
\newpage

\newpage

\section{Linearized ADMM Formulation}
Two of the numerical examples utilize variations of ADMM. This section is dedicated to a derivation of linearized ADMM used to solve problems of the form 
\begin{equation}
    \min_{x \in \bbR^n} f(Kx) + h(x) \ \ \mbox{s.t.} \ \ \|Mx-d\|\leq \delta,
    \label{eq: L-ADMM-general-problem-01}
\end{equation}
where   $K$ and $M$ are linear operators, $\delta > 0$ is a noise tolerance, and $f$ and $h$ are proximable. 
{First, we define the proximal operator for a closed, convex, and proper function by
\begin{equation}
    \prox{f}(x) \triangleq \argmin_{z\in\bbR^n} f(z) + \dfrac{1}{2}\|z-x\|^2.
\end{equation}
Letting $\delta_\sC$ be the indicator function for a closed and convex set $\sC$ with value $0$ for $x\in\sC$ and $\infty$ otherwise, the Euclidean projection $\mbox{proj}_{\sC}$ on $\sC$ is a special case of the proximal, \ie 
\begin{equation}
    \mbox{proj}_\sC(x)
    \triangleq \prox{\delta_\sC}(x)
    = \argmin_{z\in\sC} \dfrac{1}{2}\|z-x\|^2.
\end{equation}
Next} observe (\ref{eq: L-ADMM-general-problem-01}) can be rewritten as 
\begin{equation}
    \min_{x,w} f(Kx) + h(x) +  \delta_{B(d,\delta)}(w)
    \ \ \ \mbox{s.t.} \ \ \ Mx - w = 0.
    \label{eq: L-ADMM-general-problem-02}
\end{equation}
Defining the concatenation $\xi = (p,w)$, the function
\begin{equation}
    g(\xi) \triangleq f(p) + \delta_{B(d,\delta)}(w),
\end{equation}
and $S = [K; M]$, yields
\begin{equation}
    \min_{x,\xi} h(x) +  g(\xi)
    \ \ \ \mbox{s.t.} \ \ \ Sx - \xi = 0.
    \label{eq: L-ADMM-general-problem-03}
\end{equation}
Then linearized ADMM  \cite{ryu2022large} yields
{\small 
\begin{subequations}
\begin{align}
    x^{k+1} & = \prox{\beta h}\left(x^k - \beta S^\top (\nu^k + \alpha (Sx^k - \xi^k))\right) \\
    \xi^{k+1} & = \prox{\lambda g}\left(\xi^k + \lambda ( \nu^k + \alpha (Sx^{k+1}-\xi^k))\right) \\
    \nu^{k+1} & = \nu^k + \alpha (Sx^{k+1}- \xi^{k+1}),
\end{align}  
\end{subequations}} 
where $\alpha,\beta,\lambda$ are step-sizes.
Rearranging, we obtain
{\small 
\begin{subequations}
\begin{align}
    \xi^{k+1} & = \prox{\lambda g}\left(\xi^k + \lambda ( \nu^k + \alpha (Sx^{k}-\xi^k))\right) 
    \\
    \nu^{k+1} & = \nu^k + \alpha (Su^{k}- \xi^{k+1})
    \\
    x^{k+1} & = \prox{\beta h}\left(x^k - \beta S^\top (\nu^{k+1} + \alpha (Sx^k - \xi^{k+1}))\right).
\end{align}  
\end{subequations}}
Expanding $\xi^{k+1}$ reveals block-wise updates, \ie 
\begin{subequations}
\begin{align}
    p^{k+1} & = \prox{\lambda f}(p^k + \lambda (\nu_1^k + \alpha (Kx^k - p^k)) \\
    w^{k+1} & = P_{B(d,\varepsilon)}(w^k + \lambda (\nu_2^k + \alpha(Mx^k - w^k))),
\end{align}
\end{subequations}
where $\nu^k = (\nu_1^k, \nu_2^k)$ and $P_{B(d,\varepsilon_\theta)}$ is the projection onto the Euclidean ball of radius $\varepsilon_\theta$ centered about $d$. 
Writing out expanded forms gives
{\small 
\begin{subequations}
\begin{align}
    p^{k+1} & = \prox{\lambda f} \left(p^k + \lambda(\nu_1^k + \alpha(Kx^k - p^k)) \right)
    \\
    w^{k+1} & = P_{B(d,\delta)}\left(w^k + \lambda(\nu_2^k + \alpha(Mx^k - w^k)) \right)
    \\
    \nu_1^{k+1} &= \nu_1^k + \alpha(Kx^k - p^{k+1})
    \\
    \nu_2^{k+1} &= \nu_2^k + \alpha(Mx^k - w^{k+1}) 
    \\
    x^{k+1} & = \prox{\beta h}\left(x^k - \beta S^\top (\nu^{k+1} + \alpha (Sx^k - \xi^{k+1}))\right)
\end{align}  
\end{subequations}}
Expanding the update for $x^k$ reveals
\small{
\begin{subequations}
    \begin{align}
        x^{k+1}
        & = x^k - \beta S^\top (\nu^{k+1} + \alpha (Sx^k - \xi^{k+1})) 
        \\
        &= x^k - \beta \left[ \begin{array}{c} K \\ M \end{array} \right]^\top 
        \left[\hspace{-3pt}\begin{array}{c}
            \nu_1^{k+1} + \alpha (Kx^k - p^{k+1})
            \\
            \nu_2^{k+1} + \alpha (Mx^k - w^{k+1})
        \end{array}\hspace{-3pt}\right]
        \\
        &=
        x^k - \beta  K^\top\left(   \nu_1^{k+1} + \alpha (Kx^k - p^{k+1}) \right) \\
        & - \beta 
        M^\top \left( \nu_2^{k+1} + \alpha (Mx^k - w^{k+1}) \right) 
        \\
        &= x^k - \beta   K^\top \left( 2\nu_1^{k+1} - \nu_1^k  \right) \\
        & -
       \beta M^\top \left( 2\nu_2^{k+1} - \nu_2^k \right).
    \end{align}
\end{subequations}
}
The final form we implement is the tuple of update relations
{\small 
\begin{subequations}
\begin{align}
    p^{k+1} & = \prox{\lambda f} \left(p^k + \lambda(\nu_1^k + \alpha(Kx^k - p^k)) \right)
    \\
    w^{k+1} & = P_{B(d,\delta)}\left(w^k + \lambda(\nu_2^k + \alpha(Mx^k - w^k)) \right)
    \\
    \nu_1^{k+1} &= \nu_1^k + \alpha(Kx^k - p^{k+1})
    \\
    \nu_2^{k+1} &= \nu_2^k + \alpha(Mx^k - w^{k+1})
    \\
    r^{k} & =  K^\top \left( 2\nu_1^{k+1} - \nu_1^k  \right) + M^\top \left( 2\nu_2^{k+1} - \nu_2^k \right)\\
    x^{k+1} & =\prox{\beta h}\left(x^k - \beta r^k\right).
\end{align}  
\label{eq: L-ADMM-final-form}
\end{subequations}}

\section{Supplement for Implicit Dictionary} 
\label{app: dictionary}
Observe (\ref{eq: model-sparse-recovery}) is a special case of (\ref{eq: L-ADMM-general-problem-01}), taking $h =0$, $M = A$, and $f= \|\cdot\|_1$. That is, in this case, we obtain the iteration
{\small 
\begin{subequations}
\begin{align}
    p^{k+1} & = \eta_\lambda \left(p^k + \lambda(\nu_1^k + \alpha(Kx^k - p^k)) \right)
    \\
    \nu_1^{k+1} &= \nu_1^k + \alpha(Kx^k - p^{k+1})
    \\
    \nu_2^{k+1} &= \nu_2^k + \alpha(Ax^k -d) 
    \\
    r^{k} & =  K^\top \left( 2\nu_1^{k+1} - \nu_1^k  \right) + A^\top \left( 2\nu_2^{k+1} - \nu_2^k \right)\\
    x^{k+1} & = x^k - \beta r^k,
\end{align}  
\label{eq: L-ADMM-CT-appendix}
\end{subequations}}where $\eta_\lambda$ is the shrink function.

\section{Supplement for CT Reconstruction}
\label{app: CT}
Observe (\ref{eq: model-ct}) is a special case of (\ref{eq: L-ADMM-general-problem-01}), taking $h = \delta_{[0,1]^n}$, $M = A$, and $f=f_\Omega$. That is, in this case, we obtain the iteration
{\small 
\begin{subequations}
\begin{align}
    p^{k+1} & = \prox{\lambda f_\Omega} \left(p^k + \lambda(\nu_1^k + \alpha(Kx^k - p^k)) \right)
    \\
    w^{k+1} & = P_{B(d,\delta)}\left(w^k + \lambda(\nu_2^k + \alpha(Ax^k - w^k)) \right)
    \\
    \nu_1^{k+1} &= \nu_1^k + \alpha(Kx^k - p^{k+1})
    \\
    \nu_2^{k+1} &= \nu_2^k + \alpha(Ax^k - w^{k+1}) 
    \\
    r^{k} & =  K^\top \left( 2\nu_1^{k+1} - \nu_1^k  \right) + A^\top \left( 2\nu_2^{k+1} - \nu_2^k \right)\\
    x^{k+1} & = P_{[0,1]^n}\left(x^k - \beta r^k\right).
\end{align}  
\label{eq: L-ADMM-CT}
\end{subequations}}
TV minimization is obtained from (\ref{eq: model-ct})  by letting $f_\Omega$ be the $\ell_1$ norm and $K$ be a discrete differencing operator.   
For comparison to an analytic method, we use anisotropic TV minimization, \ie 
\begin{equation}
    \min_{u\in[0,1]^n} \|Du\|_1 \ \ \mbox{s.t.}\ \ \|Au-d\|\leq \varepsilon,
\end{equation}
where $\varepsilon$ is hand-tuned. The Operator Discretization Library (ODL) Python library \cite{adler2017odl} is used to compute the filtered backprojections.

\section{Supplement for Cryptoasset Trades}
\label{app: crypto}

This section is broken into three parts.
The geometric constraint sets are of particular importance to handle in a decoupled fashion, and so the first subsection is dedicated to handling CFMM constraints for batches of transactions.
Note closed form expressions exist for pairwise swaps.
The second subsection then identifies the projection operations needed.
This is followed by a derivation of a particular operator splitting used to solve the problem, giving explicit lists for updates.

\subsection{Constraint Formulation}
The dimension of the vector space for each CFMM  may differ since some exchanges might not provide access to particular cryptoassets. Consequently, we follow similarly to recent work \cite{angeris2021optimal} in using matrices $A^j \in \bbR^{n_j\times n}$ to convert global coordinates into the local coordinates of the $j$-th CFMM,
\ie 
\begin{equation}
    A^j_{k\ell} \triangleq 
    \begin{cases}
    \begin{array}{cl}
         \multirow{2}*{1} & \mbox{if token $k$ is in the $j$-th CFMM's coordinates} \\
        & \mbox{is token $\ell$ is in global coordinates}\\[5 pt]
        0 & \mbox{otherwise.}
    \end{array}
    \end{cases}
\end{equation}
Note here we use the \textit{backwards} of the referenced work, mapping global to local rather than local to global.
Let $d \in \bbR^{n\times m}$ be a matrix with the $j$-th column $d^j$ the reserve assets in the $j$-th CFMM. 
For weighted geometric CFMMs, set
\begin{equation}
    \hat{d} \triangleq  (1 + \delta)\odot d
\end{equation}
and
\begin{equation}
    \alpha_j  \triangleq    \prod_{j=1}^{n_j} (\hat{d}_i^j)^{w_i^j},
    \ \ \mbox{for all $j\in[m]$,}
\end{equation}
where $\delta_j\geq 0$ is a tolerance, $n_j$ is the number of asset types in the $j$-th CFMM, $w^j \in \bbR^{n_j}$ is a positive weighting,
and
\begin{equation}
    {\sA_j} \triangleq \left\lbrace v \in \bbR^{n_j}: v + {d}^j  \geq 0,\  \prod_{i=1}^{n_j} (v + d^j)^{w_i} \geq \alpha_j \right\rbrace.
\end{equation}
The set $\sA_j$ identifies a weighted geometric mean inequality that must hold for the $j$-th CFMM.
We include the nonnegative $\delta_j$ to account for noisy data. Choosing $\delta_j > 0$ gives a buffer for ensuring a transaction is still valid for noisy $d$ (at the cost of reducing the achievable utility $U$).

\begin{remark}
    Ideally, we would directly compute $P_{\sA_j}(x)$ in an algorithm computing optimal trades, which can be derived following an example in Beck's text. \cite{beck2017first}
    However, this projection introduces unscalable coupling since, using $\delta_j=0$ and $d=r$,
    \begin{align}
        [P_{\sA_j}(x)]_i
        & = \begin{cases}
            x_i & \mbox{if $x \in \sA_j$} \\    
             \dfrac{x_i - r_i^j +\sqrt{(x_i + r_i^j)^2+4\lambda w_i^j} }{2} & \mbox{otherwise,}
        \end{cases}
    \end{align}
    where $\lambda > 0$ is a solution to
    \begin{equation}
        \sum_{i=1}^{n_j} w_i^j \log\left(\dfrac{x_i - r_i^j+ \sqrt{(x_i+d^j_i)^2+4\lambda w_i^j}}{2}\right) = \log \alpha.
        \label{eq: lambda-projection-coupled}
    \end{equation}  
    As the number of asset types in CFMMs increase, the time of a root finding algorithm to estimate $\lambda$ also increases. Our alternative approach avoids this scaling issue. We also note JFB would technically require backpropping through the root finding scheme, but we suspect this could be avoided (by not attaching gradients during root finding) without adverse results.
\end{remark}

Upon taking logarithms, we may equivalently write
\begin{subequations}
    \begin{align}
    \sA_j 
    & = \left\lbrace v : v+d^j  \geq 0, \ \sum_{i=1}^{n_j} w^j_i \ln(v_i + d_i^j) \geq \ln(\alpha_j)\right\rbrace \\
    & = \left\lbrace v : v+d^j \geq 0, \ \left< w, \ln(v+d^j)\right> \geq \ln(\alpha_j)\right\rbrace .
    \end{align}
\end{subequations}
We decouple the constraint $\sA_j$ by defining the hyperplane 
\begin{align}
    \sH_j 
    &\triangleq \left\lbrace z : \left<w^j, z\right> = \ln(\alpha_j)\right\rbrace .
\end{align} 
and the element-wise logarithm inequality constraint set
\begin{equation}
    \sP_j \triangleq \left\lbrace (v,z) :  z \leq \ln (v + d^j), v+d^j\geq 0 \right\rbrace.
\end{equation}
These definitions yield the equivalence
\begin{equation}
    v \in \sA_j
    \ \ \iff \ \  
    \exists \ z \in \sH_j \ \mbox{s.t.} \ (v,z) \in \sP_j.
\end{equation}
This equivalence is useful since, as shown in a subsection below, $\sH_j$ and $\sP_j$ admit ``nice'' projection formulas. 
If instead the $j$-th CFMM is defined using a weighted arithmetic sum, then 
\begin{subequations}
    \begin{align}
    \sA_j 
    &=
    \left\lbrace x : x + d^j \geq 0,  \left< w, x + d^j\right> \geq  \left<w,\hat{d}^j\right> \right\rbrace \\
    & = \left\lbrace x : x + d^j \geq 0,  \left< w, x \right> \geq \left<w, \delta^j\odot d^j\right> \right\rbrace.
    \end{align}
\end{subequations} 
Next define the Cartesian product
\begin{equation}
    \sA = \sA_1\times \cdots \times \sA_m.
\end{equation}
This enables the constraints to be expressed by
\begin{equation}
    \sC_\Theta(d) = \{ (x,y) \geq 0 : A^j(\gamma_j x^j - y^j) \in \sA_j\ \forall \ j\in[m]\}.
\end{equation}
The tunable weights in $\sC_\Theta(d)$ consist of the constraint tolerances $\delta_j$.   
Let us introduce an auxiliary variable $z$ and  the block diagonal matrix $A = \mbox{diag}(A^1,\ldots,A^m)$.
Additionally, let  $\sI_1 \subset [m]$ be the subset of CFMM indices with weighted geometric product constraints and $\sI_2 \triangleq [m] - \sI_1$ the remaining indices for weighted sum constraints.
We obtain feasibility if and only if $(v,x,y,z)$ is a minimizer of the sum of indicator functions 
\begin{align}
    & \delta_{\geq 0}(x) + \delta_{\geq 0}(y) + \delta_{\sR}(v,x,y)  
    +  \sum_{j\in\sI_1} \delta_{\sP_j}(v^j,z^j) +  \delta_{\sH_j}(z^j) + \sum_{j\in\sI_2} \delta_{\sA_j}(v^j), 
\end{align} 
where
\begin{equation}
    \sR
    \triangleq \{ (v,x,y) : v = \Gamma Ax-Ay\}.
\end{equation}
This formulation of the constraints will be used in our operator splitting scheme.

\subsection{Proximal/Gradient Operations}

This section provides explicit formulas for the proximal and gradient operations needed.
First note
\begin{equation}
    P_{\geq 0}(x) = [x]_+
    \triangleq \max(x, 0),
\end{equation}
where the maximum occurs element-wise.
The projection onto a hyperplane $\sH_j$  is given by
\begin{equation}
    P_{\sH_j}(z)
    = z - \dfrac{\left<w^j,z - \ln(\hat{d}^j) \right>}{\|w^j\|^2} w^j.
\end{equation}
Similarly, if the $j$-th CFMM uses a weighted arithemtic,
\begin{equation}
    P_{\sA_j}(z)
    = z - \dfrac{[\left<w^j,z - \delta^j\odot  d^j\right>   ]_-}{\|w^j\|^2} w^j,
\end{equation}
where $[z]_- \triangleq \min(z, 0)$.
Next, the projection $P_{\sP_j}$ is defined element-wise.
The element-wise slope of $\ln(v+d)$ is $1/(v+d)$. The negative reciprocal of the slope (\ie $-(v+d)$) gives the slope of the normal line passing through the projection and the point of interest.
Letting $(\overline{v}^j,\overline{z}^j)$ be the projection of $(v^j,z^j)$ gives the point-slope relation
\begin{equation}
    \overline{z}^j - z^j = -(\overline{v}^j+d)\odot(\overline{v}^j-v^j).
    \label{eq: P-projection-point-slope}
\end{equation} 
Defining the function
\begin{equation}
    \phi(v) \triangleq {v}\odot {v} + {v}\odot ( d -v^j) - d\odot v^j + \ln({v}+d)-z^j
    \label{eq: P-projection-optimality-function}
\end{equation}
enables the relation (\ref{eq: P-projection-point-slope}) can be expressed as 
\begin{equation}
    \phi(\overline{v}^j) = 0.
\end{equation}

\renewcommand{\figurename}{Supplementary Figure}
\renewcommand{\thefigure}{S1}

\begin{figure}
    \centering
    \begin{tikzpicture}[scale=0.9]
        \draw[thick, -latex] (-1.5,0) -- (6.5,0) node[right] {$v$};
        \draw[thick] (0,-3) -- (0,4);
        \draw[thick, dashed] (-0.5,-3) -- (-0.5,4);
        \draw[thick, dashed] (2.5, -3) -- (2.5, 4);
        
        \draw[BLUE!80!black, line width=1.3] plot [domain=-0.45:6.5,samples=100] (\x,{ln(\x+0.5)}); 
        
        \draw[RED, thick, dashed] plot [domain=1.5:3.5,samples=100] (\x,{\x*\x + \x*(0.5-2.0) - 0.5*2+ln(\x+0.5) - 2.586}); 
        
        \draw[RED] (3.95,4) node[] {$\phi({v})$};
        

        \draw[BLUE!80!black] (5, 1.7047) node[above, xshift=-0.1cm, yshift=0.1cm] {$\ln(v+d)$};
        
        \draw[fill=black] (2.5, 1.086) node[below, xshift=1.2cm, yshift=-0.1cm]  {\large $P_{\sP}(v^*,z^*)$} circle(0.1) -- ($(2.5, 1.086)-0.5*(1,-3)$) circle (0.1) node[above, xshift=-0.2cm] {\large $(v^*,z^*)$};        
        
        \foreach \x in {-1,1,2,...,6}
        {
            \draw[] (\x,0.2) -- (\x, -0.2) node[below] {\x};
        }         
    \end{tikzpicture}
    \caption{Illustration for projection in $\bbR^2$ onto the set $\sP \triangleq \{(v^*,z^*) : \ln(v+d) \geq z, v+d \geq 0\}$, which is all points below the blue curve $\ln(v+d)$. Here $d=1/2$. The dashed red curve shows $\phi(v)$, the function defining the optimality condition for the projection in (\ref{eq: P-projection-optimality-function}).}
    \label{fig:my_label}
\end{figure}
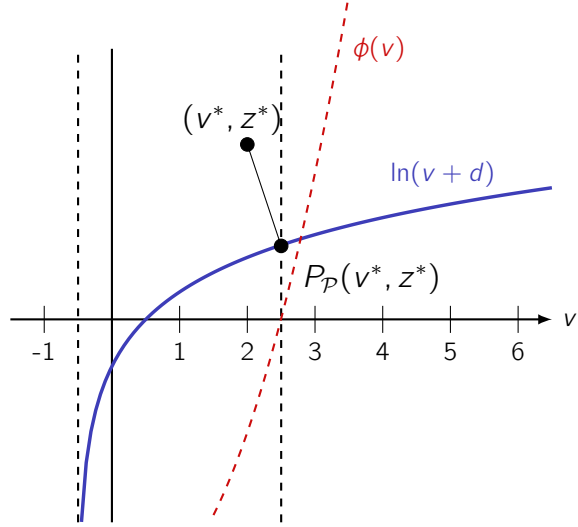

\noindent Since the above relation is element-wise and separable, each component $\overline{v}_i^j$ can be found independently (\eg via a Newton iteration).
We emphasize  solving for each $\overline{v}_i^j$ is independent of the dimension $n_j$ whereas computation costs for $\lambda$ in (\ref{eq: lambda-projection-coupled}) \textit{increase} with $n_j$.

The final projection is for the linear constraint $\sR$. The projection $P_{\sR}(v,x,y)$ is a solution to the problem
\begin{equation}
    \min_{(\overline{v},\overline{x},\overline{y})} \|\overline{v} - v\|^2
    +  \|\overline{x} - x\|^2
    +  \|\overline{y} - y\|^2
    \ \ \mbox{s.t.} \ \ 
    \overline{v} = A(\Gamma \overline{x}-\overline{y}).
\end{equation}
Let $N = [\Gamma A - A]$ and $\overline{q} = (\overline{x},\overline{y})$ so the problem becomes
\begin{equation}
    \min_{(\overline{v},\overline{q})} \|\overline{v} - v\|^2
    +  \|\overline{q} - q\|^2 
    \ \ \mbox{s.t.} \ \ 
    \overline{v} = N\overline{q}.
\end{equation}
It suffices to solve for $\overline{q}$ since the optimal $\overline{v}$ is then obtained by applying $N$.
Substituting this in yields the simpler problem
\begin{equation}
    \min_{\overline{q}} \| N\overline{q}- v\|^2
    +  \|\overline{q} - q\|^2,
\end{equation}
for which the optimality condition is
\begin{equation}
    0 = N^\top (N\overline{q}^\star - v) + \overline{q}^\star - q.
\end{equation}
Rearranging gives the formula
\begin{equation}
    \overline{q}^\star = (\mathrm{I} + N^\top N)^{-1}(q + N^\top v).
\end{equation}
Letting
\begin{equation}
    M\triangleq
    (I+N^\top N)^{-1} 
\end{equation}
and substituting in for $N^\top$ reveals
\begin{equation}
    [P_{\sR}(v,x,y)]_{(x,y)} = M\left[\begin{array}{c} x + A^\top \Gamma  v \\ y - A^\top v\end{array}\right].
\end{equation}

\noindent Lastly, we express the gradient for the utility $U_\Theta$.
Here 
\begin{align} 
    U_\Theta(x,y) 
    = \sum_{j=1}^m \left<  A^jp, A^j(y^j-x^j)\right>  
    - \dfrac{1}{2} \|W^j A^j(y^j-x^j)\|^2,
\end{align} 
where $A^j$ is used to ensure the utility only measures cryptoassets that are available on the $j$-th CFMM (\ie  converts the global coordinates of $x^j$ and $y^j$ into the local coordinates of the CFMM), and each $W^j \in \bbR^{n_j\times n_j}$ penalizes transaction sizes in the $j$-th CFMM. For each $j$,
\begin{equation}
    \nabla_{x^j} U_\Theta = -p - (W^jA^j)^\top (W^jA^j)(x^j-y^j) 
\end{equation}
and
\begin{equation}
    \nabla_{y^j} U_\Theta = -\nabla_{x^j} U_\Theta.
\end{equation}
Furthermore, $U_\Theta$ is $L$-Lipschitz with  
\begin{equation}
    L \triangleq \max_{j\in[m]} \|W^j A^j\|_2.
\end{equation}

\subsection{Operator Splitting Formulation}
Set $\xi=(v,x,y,z)$ and define the functions
\begin{subequations}
\begin{align}
    \delta_{\sM}(v,z) 
    &\triangleq   \sum_{j\in\sI_1} \delta_{\sP_j}(v^j,z^j)  + \sum_{j\in\sI_2} \delta_{\sA_j}(v^j), \\
    \delta_{\sH}(z) & \triangleq \sum_{j\in\sI_1}   \delta_{\sH_j}(z^j),
\end{align}
\end{subequations}
where $\sM$ and $\sH$ are the sets corresponding to where the indicators in their definitions are all zero.
Then define the functions
\begin{subequations}
\begin{align}
    f(\xi) & \triangleq  \delta_{\geq 0}(x,y) +\delta_{\sM}(v,z),\\
    g(\xi) & \triangleq  \delta_{\sR}(v,x,y) +  \delta_{\sH}(z)\\ 
    h(\xi) & \triangleq  -U_\Theta(x,y).
\end{align}
\end{subequations}
The problem (\ref{eq: model-crypto}) may be equivalently expressed by
\begin{equation}
    \min_{\xi} f(\xi) + g(\xi) + h(\xi),
\end{equation}
where we note $f$ and $g$ are proximable and $h$ is $L$-Lipschitz differentiable. We use Davis-Yin splitting \cite{davis2017three} and $\alpha > 0$ iterate via
\begin{subequations}
    \begin{align}
    \xi^{k+1} & = \prox{\alpha f}(\zeta^k) \\
    \psi^{k+1} & = \prox{\alpha g}(2\xi^{k+1} - \zeta^k - \alpha \nabla h(\xi^{k+1})) \\
    \zeta^{k+1} & = \zeta^k + \psi^{k+1} - \xi^{k+1}.
    \end{align}
\end{subequations}
For step size $\alpha \in (0,2/L)$, we obtain the desired convergence $(\xi_x^{k+1},\xi_y^{k+1})\rightarrow(x_d,y_d)$.

\end{document}